\pdfoutput=1
\RequirePackage{ifpdf}
\ifpdf % We are running pdfTeX in pdf mode
\documentclass[pdftex]{sigma}
\else
\documentclass{sigma}
\fi

\DeclareMathOperator{\Aut}{Aut}
\DeclareMathOperator{\Lie}{Lie}

\begin{document}

\allowdisplaybreaks

\renewcommand{\PaperNumber}{013}

\FirstPageHeading

\ShortArticleName{Principal Bundles with a~Non-Connected Structure Group}

\ArticleName{Semistability of Principal Bundles on a~K\"ahler
\\
Manifold with a~Non-Connected Structure Group}

\Author{Indranil BISWAS~$^\dag$ and Tom\'as L.~G\'OMEZ~$^\ddag$}

\AuthorNameForHeading{I.~Biswas and T.L.~G\'omez}

\Address{$^\dag$~School of Mathematics, Tata Institute of Fundamental Research,\\
\hphantom{$^\dag$}~Homi Bhabha Road, Bombay 400005, India}
\EmailD{\href{mailto:indranil@math.tifr.res.in}{indranil@math.tifr.res.in}}

\Address{$^\ddag$~Instituto de Ciencias Matem\'aticas (CSIC-UAM-UC3M-UCM),\\
\hphantom{$^\ddag$}~Nicol\'as Cabrera 15, Campus Cantoblanco UAM, 28049 Madrid, Spain}
\EmailD{\href{mailto:tomas.gomez@icmat.es}{tomas.gomez@icmat.es}}

\ArticleDates{Received October 29, 2013, in f\/inal form February 07, 2014; Published online February 12, 2014}

\vspace{-1.5mm}

\Abstract{We investigate principal $G$-bundles on a~compact K\"ahler manifold, where $G$ is a~complex
algebraic group such that the connected component of it containing the identity element is reductive.
Def\/ining (semi)stability of such bundles, it is shown that a~principal $G$-bundle $E_G$ admits an
Einstein--Hermitian connection if and only if $E_G$ is polystable.
We give an equivalent formulation of the (semi)stability condition.
A question is to compare this def\/inition with that of~[G{\'o}mez~T.L., Langer~A., Schmitt~A.H.W., Sols~I.,
  \textit{Ramanujan Math. Soc. Lect. Notes Ser.}, Vol.~10, Ramanujan Math.
  Soc., Mysore, 2010, 281--371].}


\Classification{53C07; 14F05}
\vspace{-2.5mm}

\section{Introduction}
%\label{sec1}

Let $X$ be a~compact connected K\"ahler manifold equipped with a~K\"ahler form $\omega$.
Let $G$ be a~connected complex reductive group.
The connected component, containing the identity element, of the center of $G$ will be denoted by $Z_0(G)$.
In the study of principal $G$-bundles on $X$, it is usually assumed that the group $G$ is connected.

In~\cite{Ra1}, Ramanathan gave the following def\/inition.
A principal $G$-bundle $E$ is stable (respectively, semistable) if for all proper parabolic subgroups~$P$
and all reductions of structure group $E_P\subset E|_U $ on a~big open subset $U \subset X$ (complement of
a~Zariski closed subset of codimension at least two), and for all strictly dominant characters
$\chi:P\longrightarrow \mathbb{C}^*$, the associated line bundle $E_P(\chi)$ over $U$ satisf\/ies the
inequality
\begin{gather*}
\deg(E_P(\chi))<0
\qquad
\textup{(respectively, $\deg(E_P(\chi))\leq0$)},
\end{gather*}
where the degree is calculated with respect to the K\"ahler form $\omega$.
We recall that a~strictly dominant character of $P$ is a~character of $P$ trivial on $Z_0(G)$ such that the
dual of the line bundle on $G/P$ associated to the character is ample.

A reduction $E_P\subset E\vert_U$ on a~big open set $U\subset X$ to a~parabolic subgroup $P$ is called
admissible if for every character~$\chi$ of~$P$ trivial on~$Z_0(G)$, we have $\deg(E_P(\chi)) = 0$.
A semistable principal $G$-bundle on $X$ is called polystable if there exists a~reduction $E_{L(P)} \subset
E_G$ to a~Levi factor~$L(P)$ of a~parabolic subgroup $P \subset G$ such that the following two conditions
hold:
\begin{itemize}\itemsep=0pt
\item the principal $L(P)$-bundle $E_{L(P)}$ is stable, and
\item the reduction of $E_G$ to $P$ obtained by extending the structure group of $E_{L(P)}$ to $P$ using
the inclusion $L(P) \hookrightarrow P$ is admissible.
\end{itemize}

Ramanathan used this notion in~\cite{Ra2-I,Ra2-II} to construct the moduli space of semistable principal bundles
when $\dim X=1$.
In~\cite{RS} Ramanathan and Subramanian proved that a~principal $G$-bundle admits an Einstein--Hermitian connection if and only
if it is polystable.

Behrend def\/ined stability for group schemes in~\cite{Be}.
A reductive group scheme ${\mathcal{G}}/X$ on $X$ is stable (respectively, semistable) if for all parabolic
subgroup schemes ${\mathcal{P}}/X$, we have
\begin{gather*}
\deg({\mathcal{P}}/X):=\deg({\rm ad}({\mathcal{P}}))<0
\qquad
\textup{(respectively, $\deg({\rm ad}({\mathcal{P}}))\leq0$)},
\end{gather*}
where ${\rm ad}({\mathcal{P}})$ is the associated Lie algebra bundle on $X$.
Then he def\/ines a~principal $G$-bundle on $X$ to be stable (respectively, semistable) if the associated
group scheme ${\rm Ad}(E)$ is semistable (respectively, semistable).

If $G$ is connected, then it is easy to see that Behrend's def\/inition of stability coincides with
Ramanathan's, because the normalizer $N_G(P)$ of a~parabolic subgroup $P$ in $G$ is equal to $P$ when $G$
is connected.

If $G$ is reductive but not connected, the moduli space of principal $G$-bundles has been constructed
in~\cite{GLSS} for projective varieties.
Given a~1-parameter subgroup $\lambda:\mathbb{C}^*\longrightarrow [G,G]\subset G$, a~parabolic subgroup
$P(\lambda)$ is def\/ined as follows~\cite{S}
\begin{gather}
%\label{eq:par}
P(\lambda)=\big\{g\in G:\lim_{t\to\infty}\lambda(t)g\lambda(t)^{-1}
\;
\text{exists in}
\;
G\big\}
\end{gather}
and the condition of stability is checked only with these parabolic subgroups.
The construction in~\cite{GLSS} is done using Geometric Invariant Theory, and one obtains a~stability
condition involving Hilbert polynomials.
Looking only at the leading coef\/f\/icients one obtains the associated ``slope'' stability, as usual,
which only involves degrees.
This is the stability we consider, since this is the one which is expected to correspond to the existence
of Einstein--Hermitian connections.

Here we address the problem of f\/inding a~condition for the existence of an Einstein--Hermitian connection
on a~principal $G$-bundle on $X$ when $G$ is not connected.

If $G_0$ is the connected component of identity of $G$, then $E\longrightarrow E/G_0$ is a~principal
$G_0$-bundle on $Y:=E/G_0$, and $Y$ is a~f\/inite \'etale cover of $X$.
In Section~\ref{sec2} we def\/ine the principal $G$-bundle~$E$ to be polystable when the principal
$G_0$-bundle $E\longrightarrow Y$ is polystable (cf.\ Def\/inition~\ref{def1}), and prove that this is
a~necessary and suf\/f\/icient condition for the existence of an Einstein--Hermitian connection on $E$.
In Section~\ref{sec3} we show that this def\/inition of polystability is equivalent to checking the
usual condition for those parabolic subgroups of $G$ of the form $N_G(\mathfrak{p})$, where $\mathfrak{p}$
is a~parabolic subalgebra of the Lie algebra $\mathfrak{g}$ of~$G$.
Note that $N_G({\rm Lie}(P)) = N_G(P^0)$ for any parabolic subgroup~$P$, where ${\rm Lie}(P)$ is the Lie
algebra and $P^0$ is the connected component of~$P$ containing the identity element.

We remark that, if~$G$ is connected, all parabolic subgroups are of this form because for a~connected
reductive group, $N_G(\mathfrak{p}) = N_G(P) = P$.

It is natural to ask whether our condition is equivalent to the condition in~\cite{GLSS} in terms of
\mbox{1-parameter} subgroups.
For any parabolic subgroup $P$, it is $P \subset N_G({\rm Lie}(P))$.
Hence, a~reduction of structure group to $P$ gives reduction to $N_G({\rm Lie}(P))$, and therefore if
a~principal $G$-bundle is (semi)stable in the sense of~\cite{GLSS}, then it is also (semi)stable in our
sense.

{\sloppy The implication in the other direction is not clear, since there exist examples of non-connected
groups~$G$ and 1-parameter subgroups~$\lambda$ such that~$N_G(P(\lambda)^0)$ is strictly larger
than~$P(\lambda)$.
It would be interesting to be able to compare the two def\/initions.

}

\section{Connections on principal bundles}
\label{sec2}

\subsection{Semistable and polystable bundles}

Let $G$ be a~reductive linear algebraic group def\/ined over $\mathbb C$.
We do not assume that $G$ is connected.
Let
%\begin{gather*}%\label{e1}
$G_0\subset G$
%\end{gather*}
be the connected component containing the identity element.
We note that $G_0$ is a~normal subgroup of $G$.
The quotient group
\begin{gather}
\label{e2}
\Gamma:=G/G_0
\end{gather}
parametrizes the connected components of $G$.

Let $X$ be a~compact connected K\"ahler manifold equipped with a~K\"ahler form $\omega$.
Let
\begin{gather*}
E_G \ \longrightarrow \  X
\end{gather*}
be a~holomorphic principal $G$-bundle on $X$.
Consider the quotient map
\begin{gather}
\label{e3}
E_G \stackrel{\phi}{ \ \longrightarrow \ } E_G/G_0=:Y.
\end{gather}
The natural projection
\begin{gather}
\label{e4}
f: \ Y \ \longrightarrow \  X
\end{gather}
is a~unramif\/ied Galois covering map with Galois group $\Gamma$ (def\/ined in~\eqref{e2}).
The pulled back form~$f^*\omega$ is a~K\"ahler form on $Y$.
It should be clarif\/ied that $Y$ need not be connected.

The projection $\phi$ in~\eqref{e3} makes $E_G$ a~holomorphic principal $G_0$-bundle on $Y$.
\begin{definition}
\label{def1}
The principal $G$-bundle $E_G$ on $X$ is called {\em semistable} (respectively, {\em stable}) if for each
connected component $Y'$ of $Y$, the principal $G_0$-bundle
\begin{gather*}
\phi^{-1}(Y') \ \longrightarrow \  Y'
\end{gather*}
is semistable (respectively, stable).
Similarly, $E_G$ on $X$ is called {\em polystable} if the principal $G_0$-bundle $\phi^{-1}(Y')
\longrightarrow Y'$ is polystable for every connected component $Y'$ of $Y$.
\end{definition}
\begin{lemma}
%\label{lem1}
A principal $G$-bundle $E_G$ on $X$ is semistable $($respectively, polystable$)$ if for some connected
component $Y'$ of $Y$, the principal $G_0$-bundle
\begin{gather*}
\phi^{-1}(Y') \ \longrightarrow \  Y'
\end{gather*}
is semistable $($respectively, polystable$)$.
The same criterion holds for stability.
\end{lemma}

\begin{proof}
Take two connected components $Y_1$ and $Y_2$ of $Y$.
Since the covering map $f$ in~\eqref{e4} is Galois, there is an element $g$ of the Galois group such that
the automorphism $g$ of $Y$ takes $Y_1$ to $Y_2$.
Let
\begin{gather*}
\widetilde{g}:=g\vert_{Y_1}: \ Y_1 \ \longrightarrow \  Y_2
\end{gather*}
be this isomorphism.
Let $E_1$ (respectively, $E_2$) be the restriction of the principal $G_0$-bundle $E_G \longrightarrow Y$ to
$Y_1$ (respectively, $Y_2$).
Since $f\circ \widetilde{g} = f$, and $\Lie(G) = \Lie(G_0)$, we have
\begin{gather}
\label{f1}
\widetilde{g}^*\text{ad}(E_2)=\text{ad}(E_1).
\end{gather}

On the other hand, a~principal $G_0$-bundle is semistable (respectively, polystable) if and only if the
corresponding adjoint vector bundle is semistable (respectively, polystable);
see~\cite[p.~214, Proposition 2.10]{AB} and~\cite[p.~224, Corollary 3.8]{AB}.
Therefore, from~\eqref{f1} we conclude that $E_1$ is semistable (respectively, polystable) if and only if
$E_2$ is semistable (respectively, polystable).

This isomorphism in~\eqref{f1} is compatible with the Lie algebra structure of the f\/ibers of the two
adjoint bundles.
We recall that a~principal $G_0$-bundle $F_{G_0}$ is stable if for every parabolic subalgebra bundle
$\widetilde{\mathfrak p} \subset {\rm ad}(F_{G_0})$, we have
\begin{gather*}
\text{degree}(\widetilde{\mathfrak p})<0
\end{gather*}
(see~\cite{Be}).
Therefore, if $E_1$ is stable, then $E_2$ is also stable.
\end{proof}

\subsection{Einstein--Hermitian connections}

Any two maximal compact subgroups of $G$ dif\/fer by an inner automorphism of $G$.
Fix a~maximal compact subgroup
%\begin{gather*}%\label{e5}
$K\subset G$.
%\end{gather*}
The quotient $G/K$ is a~contractible manifold, in particular, $G/K$ is connected.

Take a~holomorphic principal $G$-bundle $E_G$ over $X$.
A \textit{Hermitian structure} on $E_G$ is a~$C^\infty$ reduction of structure group of $E_G$ to $K$.
Since $G/K$ is contractible, and any $C^\infty$ f\/iber bundle with a~contractible f\/iber is trivial, it
follows immediately that $E_G$ admits Hermitian structures.

Any two connections on the principal $G$-bundle $E_G$ dif\/fer by a~smooth $1$-form with values in
$\text{ad}(E_G)$.
Two $C^\infty$ connections $\nabla_1$ and $\nabla_2$ on the principal $G$-bundle $E_G$ are called
\textit{equivalent} if $\nabla_1-\nabla_2$ is of type $(1,0)$~\cite[p.~87]{Kos}.
The complex structure on the total space of $E_G$ def\/ines an equivalence class of connections on
$E_G$~\cite[p.~87, Proposition 2]{Kos}.

Let $E_K \subset E_G$ be a~Hermitian structure.
Then there is a~unique connection $\nabla$ on the principal $K$-bundle $E_K$ such that the connection
$\widetilde{\nabla}$ on $E_G$ induced by $\nabla$ lies in the equivalence class of connections def\/ined by
the complex structure on $E_G$~\cite[pp.~191--192, Proposition 5]{At}.
This $\widetilde{\nabla}$ is called the \textit{Chern connection} corresponding to $E_K$.

Let $\text{ad}(E_K) \longrightarrow X$ be the adjoint vector bundle for the principal $K$-bundle $E_K$.
Let
\begin{gather*}
{\mathcal K}(\nabla)\in\Omega^{1,1}(\text{ad}(E_K))
\end{gather*}
be the curvature of the above connection $\nabla$ on $E_K$.
The Hermitian structure $E_K$ is called \textit{Einstein--Hermitian} if the section
\begin{gather*}
\Lambda_\omega{\mathcal K}(\nabla)\in\Omega^{0}(\text{ad}(E_K))
\end{gather*}
is given by some element of the center of the Lie algebra of $K$; here $\Lambda_\omega$ is the adjoint of
multiplication of dif\/ferential forms by the K\"ahler form $\omega$.
\begin{theorem}
\label{thm1}
A holomorphic principal $G$-bundle $E_G$ on $X$ admits an Einstein--Hermitian structure if and only if
$E_G$ is polystable.

Given a~polystable principal $G$-bundle $E_G$ over $X$, the Chern connection on $E_G$ corresponding to
a~Einstein--Hermitian structure on $E_G$ is independent of the choice of Einstein--Hermitian structure.
\end{theorem}
\begin{proof}
First assume that $E_G \longrightarrow X$ admits an Einstein--Hermitian structure.
The connection on the adjoint vector bundle $\text{ad} (E_G)$ induced by an Einstein--Hermitian connection
on $E_G$ is also Einstein--Hermitian.
Therefore, $\text{ad}(E_G)$ is polystable.
Take any connected component~$Y'$ of~$Y$.
Let $\overline{f}$ be the restriction of~$f$ to~$Y'$.
Let $E' \longrightarrow Y'$ be the principal $G_0$-bundle obtained by restricting $E_G \longrightarrow Y$
to~$Y'$.

Since $\overline{f}^*\text{ad}(E_G) = \text{ad}(E')$, and $\text{ad}(E_G)$ is polystable, we conclude that
$\text{ad}(E')$ is polystable~\cite[p.~439, Proposition 2.3]{BS}.
Hence $E'$ is polystable~\cite[p.~224, Corollary~3.8]{AB}.

To prove the converse, assume that the principal $G_0$-bundle $E_G \longrightarrow X$ is polystable.
Take a~connected component $Y'$ of $Y$.
As before, $\overline{f}$ is the restriction of $f$ to $Y'$, and
%\begin{gather*}
$E'  \longrightarrow   Y'$
%\end{gather*}
is the principal $G_0$-bundle obtained by restricting $E_G \longrightarrow Y$ to $Y'$.
The adjoint vector bund\-le~$\text{ad}(E')$ is polystable because $E'$ is polystable.
Since $\overline{f}$ is an \'etale covering map, an Einstein--Hermitian connection on $\text{ad}(E')$
produces an Einstein--Hermitian connection on the direct ima\-ge~$\overline{f}_*\text{ad}(E')$.
We note that this uses the fact that the K\"ahler form on $Y'$ is the pullback of the K\"ahler form on~$X$.

Since $\overline{f}_*\text{ad}(E')$ admits an Einstein--Hermitian connection, it follows that
$\overline{f}_*\text{ad}(E')$ is poly\-stable.
Since $\text{ad}(E_G)$ is a~direct summand of $\overline{f}_*\text{ad}(E')$, we conclude that
$\text{ad}(E_G)$ is polystable.

The Einstein--Hermitian connection on $\overline{f}_*\text{ad}(E')$ preserves $\text{ad}(E_G)$, because
$\text{ad}(E_G)$ is a~direct summand of $\overline{f}_*\text{ad}(E')$.
The connection $\nabla$ on $\text{ad}(E_G)$ obtained this way is Einstein--Hermitian.
Take the Einstein--Hermitian connection on $\text{ad}(E')$ to be one given by an Einstein--Hermitian
connection on~$E'$.
Therefore, the Einstein--Hermitian connection on~$\text{ad}(E')$ is compatible with the Lie algebra
structure of the f\/ibers of $\text{ad}(E')$.
This implies that the above connection $\nabla$ on $\text{ad}(E_G)$ is compatible with the Lie algebra
structure of the f\/ibers of $\text{ad}(E_G)$.
Therefore, $\nabla$~def\/ines a~connection on the principal $\Aut(\Lie(G))$-bundle
\begin{gather*}
E_{\Aut(\Lie(G))}=E_G\times^G\Aut(\Lie(G))
\end{gather*}
associated to $E_G$ for the homomorphism $G \longrightarrow \Aut(\Lie(G))$ given by the
adjoint action of~$G$ on~$\Lie(G)$.
This connection on~$E_{\Aut(\Lie(G))}$ given by $\nabla$ will be denoted by~$\nabla_0$.
This connection~$\nabla_0$ is Einstein--Hermitian because~$\nabla$ is so.

Def\/ine $G^{ab}:= G/[G,G]$.
The quotient homomorphism $G \longrightarrow G^{ab}$ will be denoted by $q$.
Let $G^{ab}_0\subset G^{ab}$ be the connected component containing the identity element.
Let
\begin{gather*}
\beta: \ G^{ab} \ \longrightarrow \  G^{ab}
\end{gather*}
be the homomorphism def\/ined by $z \longmapsto z^n$, where $n$ is the order of the quotient group
$G^{ab}/G^{ab}_0$.
Note that $\beta(G^{ab}_0) = G^{ab}_0$, and the homomorphism
\begin{gather*}
G^{ab}/G^{ab}_0 \ \longrightarrow \  G^{ab}/G^{ab}_0
\end{gather*}
given by $\beta$ is the trivial homomorphism.
Hence $\beta(G^{ab}) = G^{ab}_0$.
Def\/ine
\begin{gather*}
\gamma:=\beta\circ q: \ G \ \longrightarrow \  G^{ab}_0.
\end{gather*}

Since $G^{ab}_0$ is a~torus, the principal $G^{ab}_0$-bundle $E_{G^{ab}_0}$ on $X$ obtained by extending
the structure group of $E_G$ using $\gamma$ has a~unique Einstein--Hermitian connection.
We will denote this Einstein--Hermitian connection on~$E_{G^{ab}_0}$ by~$\nabla'$.

The connection $\nabla_0$ (respectively, $\nabla'$) is a~$1$-form on the total space
$E_{\Aut(\Lie(G))}$ (respectively,~$E_{G^{ab}_0}$) with values in the Lie algebra
$\Lie(\Aut(\Lie(G)))$ (respectively, $\Lie(G^{ab}_0)$).
Using the natural map $E_G \longrightarrow E_{\Aut(\Lie(G))}$ (respectively, $E_G
\longrightarrow E_{G^{ab}_0}$), the $1$-form $\nabla_0$ (respectively,~$\nabla'$) pulls back to a~$1$-form
on $E_G$ with values in $\Lie(\Aut(\Lie(G)))$ (respectively, $\Lie(G^{ab}_0)$);
this $1$-form on $E_G$ will be denoted by $\widehat\nabla$ (respectively, $\widehat{\nabla'}$).
Note that
\begin{gather*}
\Lie(G)=\Lie(\Aut(\Lie(G)))\oplus\Lie\big(G^{ab}_0\big).
\end{gather*}
Therefore, $\widehat{\nabla}+\widehat{\nabla'}$ is a~$1$-form on $E_G$ with values in the Lie algebra
$\Lie(G)$.
It is straightforward to check that this $\Lie(G)$-valued $1$-form
$\widehat{\nabla}+\widehat{\nabla'}$ def\/ines a~connection on $E_G$.
This connection on~$E_G$ will be denoted by $\widetilde\nabla$.

Fix a~point $x_0 \in X$.
Let $\Aut((E_G)_{x_0})$ denote the group of automorphisms of $(E_G)_{x_0}$ that commute with the
action of $G$ on $(E_G)_{x_0}$.
Note that $\Aut((E_G)_{x_0})$ is identif\/ied with the f\/iber of the adjoint bundle
$\text{Ad}(E_G)_{x_0}$, and it is isomorphic to~$G$.
Consider parallel translations of the f\/iber $(E_G)_{x_0}$, along loops based at $x_0$, with respect to
the above connection $\widetilde\nabla$.
These together produce a~subgroup of $\Aut((E_G)_{x_0})$.
It can be shown that this subgroup is contained in a~compact subgroup of $\Aut((E_G)_{x_0})$.
Indeed, this follows from the fact that the holonomies of both $\nabla_0$ and $\nabla'$ are compact.
Therefore, possibly taking an extension of structure group, we get a~Hermitian structure on~$E_G$.
(The above subgroup of $\Aut((E_G)_{x_0})$ is a~conjugate of a~subgroup of~$K$.) This Hermitian
structure is Einstein--Hermitian because both $\nabla_0$ and $\nabla'$ are~so.

An Einstein--Hermitian connection on the principal $G$-bundle $E_G \longrightarrow X$ pulls back to an
Einstein--Hermitian connection on the principal $G_0$-bundle $E' \longrightarrow Y'$.
Therefore, the uniqueness of the Einstein--Hermitian connection on $E_G$ follows from the uniqueness of
the Einstein--Hermitian connection on~$E'$.
To explain this, from~\cite[p.~24, Theorem~1]{RS} we know that a~stable bundle has a~unique
Einstein--Hermitian connection, and from~\cite[p.~111, Theorem~3.27]{Kob} we know that for any
decomposition of a~polystable vector bundle $F$ into a~direct sum of stable vector bundles, each direct
summand is preserved by any Einstein--Hermitian connection on~$F$.
Now apply this to the adjoint vector bundle $\text{ad}(E')$ and the principal $G_0/[G_0,G_0]$-bundle
associated to $E'$.
\end{proof}

\section{Equivalence of stability conditions}
\label{sec3}

For a~parabolic subalgebra $\mathfrak p$ of $\Lie(G)$, by $N_G(\mathfrak{p})$ we will denote the
subgroup of $G$ that preserves~$\mathfrak{p}$ by the adjoint action.
In this section, by a~parabolic subgroup of $G$ we will mean a~group of the form $N_G(\mathfrak{p})$ for
some parabolic subalgebra $\mathfrak p$.
As before, by $Z_0(G)$ we will denote the connected component of the center of $G$ containing the identity
element.
\begin{definition}
\label{adjoint}
A principal $G$-bundle is called {\it adjoint semistable} (respectively, {\it adjoint stable}) if for all
reductions to a~proper parabolic subgroup $P$, and all reductions of structure group $E_P\subset E\vert_U$
on a~big open subset $U\subset X$, and for all strictly dominant characters $\chi:P\longrightarrow
\mathbb{C}^*$, the associated line bundle $E_P(\chi)$ satisf\/ies
\begin{gather*}
\deg(E_P(\chi))<0
\qquad
(\text{respectively,~}\deg(E_P(\chi))\leq0),
\end{gather*}
where the degree is calculated with respect to the K\"ahler form $\omega$.
\end{definition}

A character of $P$ is called dominant if the restriction to $P_0:= P\bigcap G_0$ is dominant.

Recall that a~reduction $E_P\subset E\vert_U$ on a~big open set $U\subset X$ to a~parabolic subgroup $P$ is
called admissible if for every nontrivial character $\chi$ of $P$ trivial on $Z_0(G)$, we have
$\deg(E_P(\chi)) = 0$.

A semistable principal $G$-bundle on $X$ is called {\it adjoint polystable} if there exists a~reduction
$E_{L(P)}\subset E_G$ to a~Levi factor $L(P)$ of a~parabolic subgroup $P\subset G$ such that the
following two conditions hold:
\begin{itemize}\itemsep=0pt
\item the principal $L(P)$-bundle $E_{L(P)}$ is stable, and

\item the reduction of $E_G$ to $P$ obtained by extending the structure group of $E_{L(P)}$ to $P$ is
admissible.
\end{itemize}
\begin{lemma}%\label{lem:semistable}
A principal $G$-bundle $E$ on $X$ is adjoint semistable if and only if it is semistable in the sense of
Definition~{\rm \ref{def1}}.
\end{lemma}
\begin{proof}
Recall that a~principal $G$-bundle $E$ on $X$ induces a~principal $G_0$-bundle $E_0$ on \mbox{$Y {=} E/G_0$}.
By Def\/inition~\ref{def1}, a~principal $G$-bundle $E$ is semistable if and only if the restriction of
$E_0$ to a~connected component $Y'$ of $Y$ is semistable.
This is equivalent to ${\rm ad}(E_0)$ being semistable (cf.\ \cite[p.~214, Proposition~2.10]{AB}).

Note that ${\rm ad}(E_0)$ is isomorphic to $f^* {\rm ad}(E)$, which is $\Gamma$-equivariant.
If ${\rm ad}(E_0)$ is semistable, the it is also equivariantly semistable.
On the other hand, suppose that it is unstable.
Its Harder--Narasimhan f\/iltration is unique, so it will be equivariant, and hence ${\rm ad}(E_0)$ will be
equivariantly unstable.
Therefore, ${\rm ad}(E_0)$ is semistable if and only if it is equivariantly semistable.

Taking the quotient by $\Gamma$, this is equivalent to the vector bundle ${\rm ad}(E)$ on $X$ being
semistable.
We remark that the proof of~\cite[p.~214, Proposition~2.10]{AB} also works for disconnected groups if we
use Def\/inition~\ref{adjoint}.
So ${\rm ad}(E)$ is semistable if and only if $E$ is adjoint semistable.
\end{proof}
\begin{lemma}
A principal $G$-bundle $E$ on $X$ is adjoint polystable if and only if it is polystable in the sense of
Definition~{\rm \ref{def1}}.
\end{lemma}
\begin{proof}
If $E$ is polystable in the sense of Def\/inition~\ref{def1}, then from Theorem~\ref{thm1} it follows that
the adjoint vector bundle ${\rm ad}(E)$ is polystable.
Conversely, if ${\rm ad}(E)$ is polystable, then $f^*{\rm ad}(E)$ is polystable because an
Einstein--Hermitian connection on ${\rm ad}(E)$ pulls back to an Einstein--Hermitian connection on $f^*{\rm
ad}(E)$.
If $f^*{\rm ad}(E)$ is polystable, from~\cite[p.~224, Corollary 3.8]{AB} we know that $E$ is polystable in the sense of Def\/inition~\ref{def1}.

On the other hand, ${\rm ad}(E)$ is polystable if and only if $E$ is adjoint polystable; its proof is
identical to the proof of~\cite[p.~224, Corollary 3.8]{AB}.
\end{proof}

\subsection*{Acknowledgements}

We would like to thank B.~Conrad and A.~Nair for discussions.
The f\/irst-named author acknow\-led\-ges the support of the J.C.~Bose Fellowship.
The second-named author thanks the Tata Institute of Fundamental Research for the hospitality during
a~visit where part of this work was done.
This work was partly funded by the grant MTM2010-17389 and ICMAT Severo Ochoa project SEV-2011-0087 of the
Spanish Ministerio de Econom{\'\i}a y Competitividad.

\pdfbookmark[1]{References}{ref}
\LastPageEnding

\end{document}